\documentclass[12pt]{article}
\title{On definable $f$-generic groups and minimal flows in $p$-adically closed fields}
\author{Ningyuan Yao and Zhentao Zhang}

\usepackage{amsmath, amssymb, amsthm}    	
\usepackage{fullpage} 	
\usepackage{amscd}
\usepackage{hyperref}
\usepackage[all]{xy}
\usepackage{centernot}
\usepackage{color,xcolor}
\usepackage{ulem}

\DeclareMathOperator*{\forkindep}{\raise0.2ex\hbox{\ooalign{\hidewidth$\vert$\hidewidth\cr\raise-0.9ex\hbox{$\smile$}}}}

\newcommand{\GL}{\operatorname{GL}}

\newcommand{\st}{\operatorname{st}}

\newcommand{\id}{\operatorname{id}}

\newcommand{\dcl}{\operatorname{dcl}}
\newcommand{\tp}{\operatorname{tp}}

\newcommand{\dfg}{\operatorname{dfg}}

\newcommand{\img}{\operatorname{im}}

\newcommand{\cl}{\operatorname{cl}}

\newcommand{\Gen}{\operatorname{Gen}}
\newcommand{\WG}{\operatorname{WGen}}
\newcommand{\AP}{\operatorname{AP}}

\newcommand{\R}{\mathbb{R}}

\newcommand{\Q}{\mathbb{Q}}
\newcommand{\Qp}{\mathbb{Q}_p}
\newcommand{\pCF}{p\mathrm{CF}}

\newtheorem{theorem}{Theorem}[section] 
\newtheorem{lemma}[theorem]{Lemma}

\newtheorem{corollary}[theorem]{Corollary}
\newtheorem{fact}[theorem]{Fact}

\newtheorem{conjecture}[theorem]{Conjecture}

\newtheorem{proposition}[theorem]{Proposition}
\newtheorem{proposition-eh}[theorem]{Proposition(?)}
\newtheorem*{theorem-star}{Theorem}
\newtheorem*{theorem-star-A}{Theorem A}
\newtheorem*{theorem-star-B}{Theorem B}
\newtheorem*{conjecture-star}{Conjecture}
\newtheorem*{lemma-star}{Lemma}
\newtheorem*{claim-star}{Claim}

\theoremstyle{definition}
\newtheorem{definition}[theorem]{Definition}

\newtheorem{remark}[theorem]{Remark}

\newcommand{\Nn}{\mathbb{N}}

\newcommand{\ii}{{\Gen(S_C(\Qp))}}
\newcommand{\jj}{{\WG(S_H(\Qp))}}
\newcommand{\U}{\mathbb{U}}
\newcommand{\M}{\mathbb{M}}

\newcommand{\Oo}{\mathcal{O}}

\newcommand{\Ga}{\mathbb{G}_{\mathrm{a}}}
\newcommand{\Gm}{\mathbb{G}_{\mathrm{m}}}

\newcommand{\sq}{\subseteq}

\begin{document}
\maketitle
\begin{abstract}
Let $X$ be a definable group definable over a small model $M_0$. Recall that a global type $p$ on $X$   is definable $f$-generic over $M_0$ if every left translate of $p$ is definable over $M_0$. We call  $p$ strongly $f$-generic over $M_0$ if every left translate of $p$ does not fork over $M_0$.

Let $H$ be a group definable over the field $\Qp$ of $p$-adic numbers admitting global definable $f$-generic types over $\Qp$. We show that $H$ has unboundedly many global weakly generic types iff there is a global type $r$ on $H$ which is strongly $f$-generic over $\Qp$ and a $\Qp$-definable function $\theta$ such that $\theta(r)$ is finitely satisfiable in $\Qp$.

Recall that the $\mu$-type $\mu(x)$ on $H$ is the partial type consisting of the formulas over $\Qp$ which define open neighborhoods of the identity of $H$. We show that every global weakly generic type $r$ on $H$ is $\mu$-invariant: For any $\epsilon\models \mu$ and $a\models r$, we have $\epsilon\cdot a\models r$.

Let $G$  be groups definable over $\Qp$ such that $H$ is a normal subgroup of $G$ and $G/H$ is a definably compact group. Then we show that the weakly generic types on $G$ coincide with almost periodic types $G$ iff $G$ has boundedly many global weakly generic types.
\end{abstract}

\section{Introduction and Preliminaries}
In this paper we will continue the aspects of \cite{Pillay-Yao-min-flow},  \cite{Yao-Zhang-min-flow-I}, and  \cite{Yao-Zhang-min-flow-II}, where we consider the problem of whether weakly generic types coincide with almost periodic types for definably amenable NIP groups, especially groups definable in $o$-minimal structures and $p$-adically closed fields.

Let $\M$ be a monster model of a NIP theory $T$ and $(G,\cdot)$ a group definable in $\M$ with parameters from a small submodel $M\prec\M$. We call an $L_{\M}$-formula $\varphi(x)$ a $G$-formula if $\varphi(\M)\sq G$ and a partial type $\pi(x)$ a $G$-type if every formula in $\pi$ is a $G$-formula. For any $N\succ M$, let  $S_G(N)$  denote the space of all complete $G$-types over $N$. We can consider the action of $G$ on its type space $S_G(N)$.

A subset $Y$ of $G$ is (left) \emph{generic} if finitely many left $G$-translates of $Y$ cover $G$, and type $p\in S_G(M)$ is generic if $\phi(M)\cap G$ is generic for each $\phi(x)\in p$. If  $M$ is stable, then the generic types always exist. While in an unstable environment, the generic types may not exist. Newelski suggested the new notion of weakly generic types as a substitute for generic types. A  subset $Y$ of $G$ is (left) \emph{weakly generic} if there is a non-generic set $Z$ such that $Y\cup Z$ is generic, and type $p\in S_G(M)$ is weakly generic if $\phi(M)\cap G$ is weakly generic for each $\phi(x)\in p$. Recall that a \emph{Keisler measure} over $M$ on $X$, with $X$ a definable subset of $\M^n$, is a finitely additive probability measure on the Boolean algebra of $M$-definable subsets of $X$. We call  $G$   \emph{definably amenable} if there is a (left) $G$-invariant Keisler measure over $M$ on $G$.  Chernikov and Simon showed that:

\begin{fact}\cite{CS-DA-NIP}
 $G$ is definably amenable iff $G$ admits a global strongly $f$-generic type.
\end{fact}

Here, a global type $p\in S_G(\M)$ is called \emph{strongly $f$-generic over $M$} if $g\cdot p$ does not fork over $M$ for each $g\in G$. A global type on $G$ is called  \emph{strongly $f$-generic} if it is strongly $f$-generic over some small $N\prec \M$.

In fact, strong $f$-genericity can be described as following
\begin{fact}\label{fact-strongly-f-generic}\cite{CS-DA-NIP}
   Let $p\in S_G(\M)$ and $N\succ M$ be a small submodel of $\M$. Then $p$ is strongly $f$-generic over $N$ iff it is weakly generic and $N$-invariant.
\end{fact}
\begin{remark}\label{rmk-strongly-f-generic-inv}
    A global type does not fork over a small submodel $N$ iff it is $N$-invariant (see \cite{H-P-inv-measure}). So $p\in S_G(\M)$ is strongly $f$-generic over $N$ iff $g\cdot p$ is $N$-invariant for each $g\in G$.
\end{remark}

Any stable group is definably amenable since all of its global generic types are strongly $f$-generic. Among the strongly $f$-generic types, there are two extreme cases: fsg (finitely satisfiable generic) types and dfg (definable $f$-generic) types.  A global type $p\in S_G(\M)$ is \emph{fsg} or \emph{dfg} if there is a small submodel $M_0$ such that every left $G$-translate of $p$ is finitely satisfiable in $M_0$ or definable over $M_0$, respectively. We say that $G$ has fsg or dfg if $G$ admits a fsg or dfg type respectively. Of course, both fsg and dfg groups are definably amenable.

Recall that the \emph{type-definable connected component} of $G$, written $G^{00}$, is the smallest type-definable subgroup of $G$ of bounded index, which always exists as $G$ is an NIP group \cite{HPP-$NIP$-JAMS}. The \emph{definable connected component} of $G$, written $G^{0}$, is the intersection of all definable subgroups of $G$ of finite index. Clearly, $G^{00}\leq G^0$ are normal subgroups of $G$.  If $G$ is definable over $\Qp$, then $ G^0=G^{00}$ (see \cite{J-Y-connected-component}).
\begin{fact}\cite{CS-DA-NIP}
    Let $p\in S_G(\M)$ and $\mathrm{stab}(p):=\{g\in G|\ g\cdot p=p\}$ the stabilizer of $p$. If $G$ is definably amenable, then $p$ is weakly generic iff $ G^{00}=\mathrm{stab}(p)$
    iff $p$ has a bounded $G$-orbit.
\end{fact}

Let $M$ be the field $\Qp$ of $p$-adic numbers or an $o$-minimal expansion of a real closed field, then $G$ has dfg \cite{Conversano-Pillay, P-Y-dfg groups} iff $G$ is totally non-compact and $G$ has fsg iff $G$ is definably compact \cite{HPP-$NIP$-JAMS, O-P, Johnson-fsg}.

In \cite{Yao-Zhang-min-flow-II}, we studied the badness and $\mu$-invariance of the commutative dfg groups definable over $\Qp$. Let $M=\Qp$ and $G$ a  group definable over $\Qp$. Recall from \cite{Yao-Zhang-min-flow-I} that
 a dfg $G$ is called \emph{bad} if there is a global type $p\in S_G(\M)$ strongly $f$-generic over $\Qp$ and a $\Qp$-definable map $\theta$ such that $\theta(p)$ is finitely satisfiable but not realized in $\Qp$. This definition is motivated by the following fact:
 \begin{fact}\cite{Yao-Zhang-min-flow-I}\label{fact-heir-neq-coheir}
Let $e\in \M\backslash \Qp$  and $p\in S_1(\Qp)$ be a non-algebraic type. Suppose that $p_1$ is the unique heir  $p$ over $\Qp,e$, then $p_1$ is not finitely satisfiable in $\Qp$.
\end{fact}

 We showed in \cite{Yao-Zhang-min-flow-II} that:
\begin{fact}\label{fact-bad-stationary-commu}
    Let $G$ be a commutative dfg group definable over $\Qp$, then $G$ is bad iff $G$ has unboundedly many global weakly generic types.
\end{fact}

Note that any group $G$ definable over $\Qp$ has a definable manifold structure over $\Qp$ such that the group operation is analytic (see Lemma 3.8 of \cite{Pillay-On fields definable in Qp}).  Moreover, $G$ has a family of $\Qp$-definable open compact subgroups which forms a local base of the identity $\id_G$ of $G$.   Let $\mu_G(x)$ be the partial type of all formulas over $\Q_p$ defining open compact subgroups of $G$. We call $\mu_G(x)$ the \emph{infinitesimal type} of $G$ (over $\Qp$). Then $\mu_G(\M)$ is a type-definable subgroup over $\Qp$. Clearly, $\mu_G(\M)\subseteq  G^0$. If $G$ is definably compact, letting $\st: G\to G(\Qp)$ be the standard part map, then $\mu_G(\M)=G^0=\st^{-1}(\id_G)$ (see \cite{O-P}).

Let $N\succ\Qp$ and $p\in S_G(N)$, we say that $p\in S_G(N)$ is \emph{$\mu$-invariant} if $\tp(a\cdot b/N)=p$ for any $a\models \mu_G$ and $b\models p$. We call $G$ \emph{$\mu$-invariant} if every global  weakly generic type on $G$ is $\mu$-invariant. Note that there is no $\mu$-invariant type when $G$ has fsg. 
On the other side, for dfg groups we showed in \cite{Yao-Zhang-min-flow-II} that:
\begin{fact}\label{fact-mu-inv-commu}
    Let $G$ be a commutative dfg group definable over $\Qp$, then $G$ is $\mu$-invariant.
\end{fact}

In this paper, we generalize Fact \ref{fact-bad-stationary-commu} and Fact \ref{fact-mu-inv-commu}  to the non-commutative case:
\begin{theorem-star-A}
  Let $G$ be a dfg group definable over $\Qp$, then
  \begin{enumerate}
      \item [(i)] $G$ is bad iff $G$ has unboundedly many global weakly generic types.
      \item [(ii)]  $G$ is $\mu$-invariant.
  \end{enumerate}
\end{theorem-star-A}

We now assume that $M$ is an arbitrary structure and consider $S_G(M)$ as a $G$-flow. A subset $X\sq S_G(M)$ is a \emph{subflow} if $X$ is closed and invariant under the group action, i.e.  $G\cdot X=X$. A type $p\in S_G(M)$ is called \emph{almost periodic} if $p$ is contained in some minimal subflow, or equivalently, the closure $\cl(G\cdot p)$ of the $G$-orbit of $p$ is a minimal subflow. Let $\Gen(S_G(M))$, $\WG(S_G(M))$, and $\AP(S_G(M))$ be the space of generic types, weakly generic types, and almost periodic types, respectively. Newelski showed in \cite{Newelski-I} that
\begin{fact}\label{fact-newelski}
  \begin{enumerate}
    \item [(i)] $\AP(S_G(M))=\cl(\WG(S_G(M)))$.
    \item [(ii)] If $\Gen(S_G(M))\neq \emptyset$, then $S_G(M)$ contains a unique minimal subflow and
    \[
    \Gen(S_G(M))=\WG(S_G(M))=\AP(S_G(M)).
    \]
\end{enumerate}
\end{fact}


It is easy to see that $\Gen(S_G(M))=\WG(S_G(M))=\AP(S_G(M))$ when  $M$ is stable. Newelski raised a question to find an $o$-minimal or even NIP example where $\WG(S_G(M))$ and $\AP(S_G(M))$ differ. Later, the question was reformulated by Chernikov and Simon. They asked in \cite{CS-DA-NIP} whether or not $\WG(S_G(M))=\AP(S_G(M))$ if $G$ is a definably amenable NIP group.

When $M$ is an $o$-minimal expansion of a real closed field, then $G$ is definably amenable iff $G$ has a (definable) short exact sequence
\[
1\to H\to G\to C\to 1,
\]
where $H$ has dfg and $C$ has fsg \cite{Conversano-Pillay}. 
Pillay and the first author proved in \cite{Pillay-Yao-min-flow} that weakly generics coincide with almost periodics when the $H$ has dimension one. They also gave the example $G=(\R^2,+)\times {\mathbb S}^1$ for $\WG(S_G(M))\neq \AP(S_G(M))$. They conjectured in \cite{Pillay-Yao-min-flow} that
\begin{conjecture}\label{Conj-DA}
   Any definably amenable NIP group is an extension of an fsg group by a dfg group.
\end{conjecture}
The recent work of Jonhson and the first author shows that any commutative group definable in a $p$-adically closed field is an extension of an fsg group by a dfg group \cite{J-Y-abelian},  providing positive evidence for Conjecture \ref{Conj-DA}.

In our previous paper \cite{Yao-Zhang-min-flow-I}, we deal with the case where $M$ is $\Qp$ or an $o$-minimal expansion of a real closed field. With the assumption that  $G$ has a definable short exact sequence
\[
1\to H\to G\to C\to 1,
\]
where $H$ has dfg and $C$ has fsg, we showed that:
\begin{fact}\cite{Yao-Zhang-min-flow-I}\label{fact-G-is-stationary-H-is-stationary} Let $G$ be as above. Then
   \begin{enumerate}
    \item [(i)] If $G$ has boundedly many global weakly generic types, then $\WG(S_G(M))=\AP(S_G(M))$.
    \item [(ii)] $G$ has boundedly many global weakly generic types iff $H$ has boundedly many global weakly generic types.
\end{enumerate}
\end{fact}

When $G$ is a commutative group definable over $\Qp$, applying Fact \ref{fact-bad-stationary-commu}, Fact \ref{fact-mu-inv-commu}, and Fact \ref{fact-G-is-stationary-H-is-stationary}, we showed in \cite{Yao-Zhang-min-flow-II} that:
\begin{fact}\label{fact-commu-WG=AP-equiv}
The following are equivalent:
\begin{enumerate}
\item [(i)] $G$ has boundedly many global weakly generic types.
\item [(ii)] $\WG(S_G(\Qp))=\AP(S_G(\Qp))$.
\item [(iii)] $\WG(S_G(\M))=\AP(S_G(\M))$.
\end{enumerate}
\end{fact}
In this paper, applying Theorem A, we generalize Fact \ref{fact-commu-WG=AP-equiv} to the non-commutative case:
\begin{theorem-star-B}
Let $G$ be a group definable over $\Qp$.  If $G$ is a definable extension of a fsg group $C$ by a dfg group $H$, where $C$ is infinite, then the following are equivalent:
\begin{enumerate}
\item [(i)] $G$ has boundedly many global weakly generic types.
\item [(ii)] $\WG(S_G(\Qp))=\AP(S_G(\Qp))$.
\item [(iii)] $\WG(S_G(\M))=\AP(S_G(\M))$.
\end{enumerate}
\end{theorem-star-B}

\subsection{ Notation and conventions}
We use $\M$ to denote the monster model of $\pCF$, the complete theory of $\Qp$ in the language of rings. A subset $A$ is called \emph{small} or \emph{bounded} if $|A| < |\M|$. We use $M$, $N$, $M_1$, $N_1$,... to denote ``small'' elementary elementary extensions of $\Qp$. We usually write tuples as $a, b, x, y...$ rather than  $\bar a,\bar b,\bar x,\bar y$. By a ``global type'', we mean a complete type over $\M$.

When we speak of a set definable {\em in} $M$ or $\M$, we mean the obvious thing.  When we speak of a set $X$ definable {\em over} $M$, we typically mean a set definable in $\M$ defined with parameters from $M$. By a definable object, we mean a definable object in the monster model $\M$. We sometimes use $X(x)$ to denote the formula which defines $X$. In this case, $X(M)$ denotes the set definable in $M$ by the same formulas defining $X$ in $\M$.



\subsection{Outline}
In Section 2, we review some notions of $p$-adic algebraic groups, especially the trigonalizable algebraic group, and prove Theorem A. In Section 3, we review some notions and previous results of definable topological dynamics and prove Theorem B.

\section{Badness, stationarity, and $\mu$-invariance of dfg groups}

Let   $F$ be an algebraically closed field containing $\Qp$. Let $A$ be an algebraic group definable in $F$  over $\Qp$, namely, the variety structure as well as the group structure are given by data (polynomial equations, transition maps, morphisms) over $\Qp$ (See \cite{Pillay-ACF}). We call  $A(\Qp)$ of $\Qp$-points of the algebraic group $A$ a  {\em $p$-adic algebraic group}.  Of course, $A(\Qp)$ will be also a definable group in the structure  $(\Qp,+,\times,0,1)$, but essentially just quantifier-free definable in the ring language.  We call $A(\Qp)$ is connected if $A$ is connected as an algebraic group.

 Let $G\sq \M^n$ be a group definable over $\Qp$. We call $G$ an \emph{algebraic group over $\Qp$} if $G(\Qp)$ is a  $p$-adic algebraic group. Similarly, $G$ is connected (resp. linear, unipotent, split,...) if $G(\Qp)$ is connected (resp. linear, unipotent, split,...). We assume that all algebraic groups are connected in this paper.

Recall from \cite{Milne-book} that a \emph{trigonalizable algebraic group $H$  over $\Qp$} is a linear solvable algebraic group over $\Qp$ which is also split over $\Qp$. Equivalently, $H\leq \GL(n,\M)$  has a decompostion  $U\rtimes T$, where $U$ is unipotent   and $T$ is a torus, and both of them split over $\Qp$. Let $\Gm$ denote the multiplicative group  $(\M^*,\times)$, then  $T$ is  isomorphic to $\Gm^k$ via a (quantifier-free) $\Qp$-definable map for some $k\in\Nn$. There exists a matrix $g$ in $\GL(n,\Qp)$ such that
$gUg^{-1}$ is contained in the group $\U_n\leq \GL(n,\M)$ of upper unitriangular matrices. It follows that there is a central series
\[
U=U_0\geq U_1\geq\cdots\geq U_n=\{\id_U\}
\]
in $U$ such that $U_i/U_{i+1}\cong \Ga$ for $i=0,...,n-1$, where $\Ga=(\M,+)$ is the   additive group. If $U$ is commutative, then $U$ is definably isomorphic to a product of copies of $\Ga$'s over $\Qp$ (see section 16 of \cite{Milne-book} or section 2.1.8 of \cite{P-R-AG-book} for details).  Note that $\Ga^0=\Ga$ and $\Gm^0=\bigcap_{n\in \Nn^{>0}}P_n(\Gm)$ and $P_n(\Gm)$ is the set of $n$-th powers in $\Gm$ (see Proposition 2.3 and 2.4 of \cite{PPY-sl2}).

Let  $X $ and $A $ be groups definable over $\Q_p$. When we say that ``$X$ is \emph{virtually} $A$'', we mean that ``there are a finite index $\Q_p$-definable subgroup $X_0$ of $X$ and a $\Q_p$-definable morphism $f: X_0\rightarrow A$ such that both $\ker(f)$ and $A/\img(f)$ are finite''.
\begin{fact}\label{fact-dfg-groups}\cite{P-Y-dfg groups}
 Let $X$ be a group definable over $\Qp$. Then
 \begin{enumerate}
     \item [(i)] $X$ has $\dfg$ iff $X$ is virtually  a   trigonalizable algebraic group over $\Qp$. \cite{P-Y-dfg groups}
     \item [(ii)] If $X$ has dfg, then  $\WG(S_X(N))=\AP(S_X(N))$ for any $N\succ \Qp$ (see \cite{P-Y-dfg groups,Yao-tri}).

 \end{enumerate}
\end{fact}

So all dfg groups are virtually trigonalizable algebraic groups. Recall from \cite{Yao-Zhang-min-flow-I} that a definable group $G$ over $\Qp$ is \emph{stationary} if every weakly generic type $p\in S_G(\Qp)$ has a unique global weakly generic extension. Note that every fsg group over $\Qp$ is stationary \cite{O-P}.
\begin{fact}\cite{Yao-Zhang-min-flow-I}\label{fact-stationary}
Let $X$ be a dfg group definable over $\Qp$, then the following are equivalent:
\begin{enumerate}
    \item [(i)] $X$ is stationary.
    \item [(ii)] $X$ has boundedly many global weakly generic types.
    \item [(iii)] Every global weakly generic type on $X$ is definable over $\Qp$.
\end{enumerate}
\end{fact}

Recall  that A dfg group $X$ over $\Qp$ is bad if there is a global strongly $f$-generic type $\tp(h/\M)\in S_X(\M)$ and
 a $\Qp$-definable map $\theta$ such that $\tp(\theta(h)/\M)$ is finitely satiafiable in $\Qp$.  It is easy to see from Fact \ref{fact-stationary} that a stationary dfg group is not bad.

\begin{fact}\cite{Yao-Zhang-min-flow-II}\label{fact-virtually-sationary-bad}
Let $G_1$ and $G_2$ be dfg groups over $\Qp$. If $G_1$ is virtually $G_2$, then
\begin{enumerate}
\item [(i)] $G_1$ is stationary iff $G_2$ is.
\item [(ii)] $G_1$ is bad iff $G_2$ is.
\end{enumerate}
\end{fact}

According to Fact \ref{fact-dfg-groups} and Fact \ref{fact-virtually-sationary-bad}, if one wants to study the stationarity and badness of dfg groups over $\Qp$, it suffices to study the stationarity and badness of trigonalizable algebraic groups over $\Qp$.

In this section, $H\leq \GL(n,\M)$ will denote a trigonalizable algebraic group over $\Qp$. We may assume that $H$   has a decomposition $U\rtimes T$ with $U$ is an algebraic subgroup of $\U_n$ and $T=\Gm^k$ for some $k\in \Nn$.  Clearly, $H$ is commutative iff $H=U\times T$ and $U\cong \Ga^s$ for some $s\in\Nn$.
It is well-known that if $\dim(U)=2$ then $U$ is commutative, so it is definably isomorphic to $\Ga\times \Ga$.

\begin{fact}\cite{Yao-tri}
Let $U$ and $T$ be as above. Then
\begin{enumerate}
    \item [(i)] $U^{0}=U$, and $T^{0}={(\Gm^0)}^n$.
    \item [(ii)] $H^{0}=U\rtimes T^0$.
\end{enumerate}
\end{fact}

\begin{lemma}\label{lemma-uni-semi-pro}
  $U$ is definably isomorphic to  $U_1\rtimes \Ga$  where $U_1$  is a unipotent algebraic group over $\Qp$ of dimension $\dim(U)-1$.
\end{lemma}
\begin{proof}
Let $U_1\leq U$ be the normal subgroup such that $U/U_1\cong \Ga$. We can regard the Lie algebra  $L(U)$  of $U$  as a subspace of the Lie algebra $L(\U_n)$ of $\U_n$, which is the space of upper triangular matrices with $0$'s along the diagonal. The exponential map $\exp: L(U)\to U$ is given by
\[
\exp: g\mapsto I+g+\frac{g^2}{2!}+...+\frac{g^{n-1}}{(n-1)!},
\]
where $I$ is the identity of $U$. Clearly, $\exp$   is definable in the language of rings. Take a vector $v$ in the Lie algebra of $U$ which is not in the Lie algebra of $U_1$. Let $E$ be the image of the exponential map $\{\exp(tv) : t \in \M\}$. Then $E\cong \Ga$  and $E\not\subseteq U_1$. Now $E\cap U_1$ is a proper unipotent algebraic subgroup of $E$, so is trivial (see Exercise 14-3 of \cite{Milne-book}). Let $U'=EU_1$, then $U'$ is a unipotent algebraic subgroup of $U$ and $\dim(U_1)<\dim(U')$. We conclude that $U=U'$. So $U$ is a semidirect product of $U_1$ and $E$.
\end{proof}

 Now $H$ is a semidirect product of $U$ and $T$. A similar argument as Corollary 2.17 in \cite{Yao-Zhang-min-flow-II} shows that if either $U$ or $T$ is bad, then $H$ is bad. When $H$ is commutative, we showed in \cite{Yao-Zhang-min-flow-II} that
 \begin{fact}\label{fact-commu-stationary-bad}
 If $H$ is commutative, then
 \begin{enumerate}
  \item [(i)] $H$ is bad iff $\dim(U)\leq 1$.
     \item [(ii)] $H$ is stationary iff $H$ is not bad.
 \end{enumerate}
 \end{fact}

\begin{fact}\label{fact-gen-Gm-n}\cite{Yao-Zhang-min-flow-II}
$T=\Gm^k$ is stationary for each $k\in \Nn$. Moreover, $\tp(a_1,\dots,a_k/\M)$ is a weakly generic type on $T$ if and only if each $\tp(a_i/\M,a_j:i<j)$ is $\Gm^{0}$-invariant ($\Gm^0=\Gm^0(\M)$).
\end{fact}

Note that Fact \ref{fact-commu-stationary-bad} (i) is NOT true when $H$ is not commutative  (see Lemma 4.11 of \cite{Yao-Zhang-min-flow-I}). We now generalize Fact \ref{fact-commu-stationary-bad} (ii) to the non-commutative case.

 \begin{lemma}\label{lemma-dim(U)>1-bad}
  If $\dim(U)>1$, then $U$ is bad, so $H$ is not stationary.
 \end{lemma}
\begin{proof}
  If $\dim(U)= 2$, then $U$ is commutative, so $U$ is bad by Fact \ref{fact-commu-stationary-bad}. So we only need to consider the case $\dim(U)\geq 3$. Induction on the dimension on $U$. By Lemma \ref{lemma-uni-semi-pro}, $U=U_1\rtimes \Ga$ with $2\leq \dim(U_1)=\dim(U)-1$ , we see that $U_1$ is bad by the induction hypothesis. So $U$ is bad as required.
\end{proof}

We now consider the case where $\dim(U)=1$. We may assume that $U=\Ga$ and $T=\Gm^k$. The action $\rho: T\times \Ga\to \Ga$  given by $(t,u)\mapsto u^t$ is an $\Qp$-definable map.

\begin{lemma}\label{lemma-Qp-automorphism}
Let $f:\Ga(\Qp)\to\Ga(\Qp)$ be a $\Qp$-definable automorphism. Then  $f$ is linear. In fact, $f(a)=f(1)a$ for all $a\in\Ga(\Qp)$.
\end{lemma}
\begin{proof}
Clearly, $f$ is continuous since it is a definable morphism of definable groups.  As $\Ga(\Qp)$ is torsion-free and divisible, we see that $f$ is a $\Q$-linear map. So $f(q)=qf(1)$ for each $q\in \Q$. Since $\Q$ is dense in $\Ga(\Qp)$ and $f$ is continuous, we see that $f(a)=af(1)$ for all $a\in \Ga(\Qp)$. This completes the proof.
\end{proof}

\begin{lemma}\label{lemma-Qp-action}
Let $B$ be a group definable over $\Qp$ and $\eta: B\times \Ga\to \Ga$ a $\Qp$-definable action of $B$ on $\Ga$ such that $\eta(b,-):\Ga\to\Ga$ is an automorphsim of $\Ga$ for each $b\in B$. Then there is a $\Qp$-definable morphism $\rho: B\to \Gm$ such that $\eta(b,a)=\rho(b)a$ for each $b\in B$ and $a\in \Ga$.
\end{lemma}
\begin{proof}
For each $b\in B(\Qp)$, $\eta(b,x)$ is an automorphsim of $\Ga$.   By Lemma \ref{lemma-Qp-automorphism}, we have
\[
\Qp\models \forall b\in B(\Qp) \exists ! u\in \Gm(\Qp)\forall x\in \Ga(\Qp)(\theta(h,x)=ux),
\]
which implies that there is a $\Qp$-definable map  $\rho: B(\Qp)\to \Gm(\Qp)$ such that $\eta(b,a)=\rho(b)a$ for each $b\in B(\Qp)$ and $a\in \Ga(\Qp)$. Now for any $b_1,b_2\in B(\Qp)$ and $a\in\Ga(\Qp)$,
\[
\eta(b_1,\eta(b_2,a))=\rho(b_1)\rho(b_2)a \ \ \text{and}\ \ \eta(b_1b_2,a))=\rho(b_1b_2)a .
\]
Since $\eta(b_1,\eta(b_2,a))=\eta(b_1b_2,a))$, we see that $\rho(b_1b_2)=\rho(b_1)\rho(b_2)$. So $\rho:B(\Qp)\to\Gm(\Qp)$ is a group morphism. Working in the monster model, we see that   $\rho: B\to \Gm$ is a $\Qp$-definable group morphism such that $\eta(b,a)=\rho(b)a$ for all $b\in B$ and $a\in\Ga$.
\end{proof}

Suppose that $H=\Ga\rtimes T$ with $T=\Gm^k$. By Lemma \ref{lemma-Qp-action}, there is a definable morphism $\rho: T\to \Gm$ such that $u^t=\rho(t)u$ for all $u\in\Ga$ and $t\in T$. We will denote $H$ by \emph{$\Ga\rtimes_\rho T$} in this case. Note that if $\tp(t/\M)\in S_T(\M)$ is definable over $\Qp$, then $\tp(\rho(t)/\M)$ is also definable over $\Qp$. So $\rho(T)\leq \Gm$ is either an infinite dfg group over $\Qp$ or a finite subgroup of $\Gm(\Qp)$.

Recall the $\pCF$ is a distal theory (see \cite{Sim-distal, Sim-fin} for details). By Lemma 2.16 in \cite{Sim-distal}, we have the following fact:
\begin{fact}\label{fact-distal}
Let $N\prec \M$ a small submodel, $p(x)\in S(\M)$ definable over $N$, and $q(y)\in S(\M)$ finitely satisfiable in $N$. Then $p(x)$ and $q(y)$ are \emph{orthogonal}. Namely, $p(x)\cup q(y)$ implies a complete global type. In fact, if $a\models p$ and $b\models q$, then $\tp(a/\M,b)$ is the unique heir of $\tp(a/N)$ and $\tp(b/\M,a)$ is finitely satisfiable in $N$.
\end{fact}

\begin{lemma}\label{lemma-dim(U)=1-stationary}
  Let $H=\Ga\rtimes_\rho T$. If $\rho(T)$ is finite, then $H$ is stationary.
\end{lemma}
\begin{proof}
If $\rho(T)$ is finite, then $T_0=\ker(\rho)$ is a finite  index subgroup of $T$. So $\Ga\rtimes_\rho T_0=\Ga\times T_0$  is a finite index subgroup of $H$, so $H$ is virtually $\Ga\times \Gm^k$, thus is stationary by Fact\ref{fact-virtually-sationary-bad} and Fact \ref{fact-commu-stationary-bad}.
\end{proof}

\begin{lemma}\label{lemma-dim(U)=1-bad}
  Let $H=\Ga\rtimes_\rho T$. If $\rho(T)$ is infinite, then $H$ is bad.
\end{lemma}
\begin{proof}
Since $\rho(T)$ is an infinite dfg subgroup of $\Gm$, it has finite index in $\Gm$ (see \cite{johnson-definable type pCF}). Let $q=\tp(t^*/\M)$ be a global  weakly generic type on $T$, then $q$ is definable over $\Qp$ since $T$ is stationary. Clearly, $\tp(\rho(t^*)/\M)$ is a weakly generic type on $\Gm$. Let $\Gamma_\M$ be the value group of $\M$ and $v:\M\to\Gamma_\M\cup\{\infty\}$ the valuation map. Then we have that either $v(\rho(t^*))<\Gamma_\M$ or $v(\rho(t^*))>\Gamma_\M$ (see Proposition 2.4 of \cite{PPY-sl2}). Replacing $t^*$ by ${t^*}^{-1}$ if necessary, we may assume that  $v(\rho(t^*))<\Gamma_\M$.

 Let $\Oo$ be the valuation ring of $\M$ and $r=\tp(\epsilon/\M)$ be is a global fsg type on $(\Oo,+)$ such that $r\models \Oo^0$. Note that $\Oo^0$ is defined by $\mu_{\Ga}(x)$, the infinitesimal type of $\Ga$ over $\Qp$. Let $u^*=\epsilon\rho(t^*)$, then  $(u^*,t^*)$ and $(\epsilon,t^*)$ are interdefinable over $\Qp$. Since both $q$ and $r$ are orthogonal and both of them are $\Qp$-invariant, we see that   $\tp(\epsilon,t^*/\M)$ is also $\Qp$-invariant, and hence, $p=\tp(u^*,t^*/\M)$ is $\Qp$-invariant. If $p$ is a weakly generic type on $H$, then $p$ is strongly $f$-generic over $\Qp$ (see \cite{CS-DA-NIP}), and it witnesses the badness of $H$ since $\tp(\frac{u^*}{\rho(t^*)}/\M)$ is finitely satisfiable in $\Qp$.  We now show that $p$ is weakly generic. It suffices to show that $p$ is $H^0=\Ga\rtimes T^0$-invariant.

 Let $(u,t)\in H^0$, then
\[
(u,t)\cdot p=\tp(u+\rho(t)\epsilon\rho(t^*),tt^*/\M)=\tp(u+\epsilon\rho(tt^*),tt^*/\M).
\]
To see
\[
\tp(u+\epsilon\rho(tt^*),tt^*/\M)=\tp( \epsilon\rho(t^*), t^*/\M)=p,
\]
it suffices to show that
\[
\tp(\frac{u+\epsilon\rho(tt^*)}{\rho(tt^*)},tt^*/\M)=\tp( \epsilon, t^*/\M).
\]
Applying the orthogonality, it suffices to show that
\[
\tp(\frac{u+\epsilon\rho(tt^*)}{\rho(tt^*)}/\M)=\tp( \epsilon/\M).
\]
Since $\rho(T^{0})=\rho(T)^0\subseteq \Gm^0$ (see Lemma 2.1 of \cite{J-Y-abelian})),
\[
\tp(\rho(tt^*)/\M)=\tp(\rho(t)\rho(t^*)/\M)=\tp(\rho(t^*)/\M).\]
We conclude that $v(\rho(tt^*))<\Gamma_\M$ and thus  $tp(\frac{u}{\rho(tt^*)}/\M)$ is either $0$ (if $u=0$) or infinitesimally close to $0$ over $\M$. Applying the orthogonality again, we see that $\tp(\epsilon/\M,\frac{u}{\rho(tt^*)})$ is a fsg type on $\Oo$ over $\dcl(\M,\frac{u}{\rho(tt^*)})$. Sicne $\frac{u}{\rho(tt^*)}$ realizes $\mu_{\Ga}(x)$, we have that
\[
\tp(\frac{u}{\rho(tt^*)}+\epsilon/\M,\frac{u}{\rho(tt^*)})=\tp(\epsilon/\M,\frac{u}{\rho(tt^*)}).
\]
So
\[
\tp(\frac{u+\epsilon\rho(tt^*)}{\rho(tt^*)}/\M)=\tp( \epsilon/\M)
\]
as required.
\end{proof}

\begin{proposition}\label{prop-stationary iff not bad}
$H$ is stationary iff $H$ is not bad.
\end{proposition}
\begin{proof}
The direction of $\Rightarrow$ is obvious.  Conversely, suppose that $H=U\rtimes T$ and $H$ is not bad,  then $\dim(U)= 1$ by Lemma \ref{lemma-dim(U)>1-bad}.  We may assume that  $H$ is of the form $\Ga\rtimes_\rho T$. By Lemma  \ref{lemma-dim(U)=1-bad}, $H$ is not bad implies $\rho(T)$ is finite, which implies that $H$ is stationary by Lemma \ref{lemma-dim(U)=1-stationary}.
\end{proof}

We now consider the $\mu$-invariance of $H$. We recall some notions from \cite{P-S: Top gup}. Let $G$ be a group definable over $\Qp$.  If $\varphi(x)$ and $\psi(x)$ are $G$-formulas, then by $\varphi\cdot \psi$ we  denote the $G$-formula
        \[
        (\varphi\cdot \psi)(x)=\exists u\exists v (\varphi(u)\wedge\psi(v)\wedge x=uv).
        \]
If $p(x)$ and $r(x)$ are (partial) $G$-types, then
        \[
        (p\cdot r)=\{(\varphi\cdot \psi)(x)|\ p\vdash \varphi(x),\ r\vdash \psi(x)\}
        \]

 Let $N$ be an elementary extension of $\Q_p$. We say that $p\in S_G(N)$ is (left)  \emph{$\mu$-invariant} if $\mu_G \cdot p=p$.


We now fix $N\succ \Qp$ as an arbitrary  (small) submodel of $\M$ which is also sufficiently saturated.
\begin{definition}
 We say that $G$ has the \emph{$\sharp$-property} if for any $N^*\succ N$, each $H^{0}(N)$-invariant type $p\in S_G(N^*)$ is $\mu$-invariant.
\end{definition}

\begin{fact}\label{fact-sharp-virtually-mu-inv}[ Lemma 3.6 of \cite{Yao-Zhang-min-flow-II}]
Suppose that $A$ and $B$ be defianbly  amenable groups definable over $\Qp$.  Then we have:
\begin{enumerate}
    \item [(i)] If $A$ is virtually  $B$, then $A$  has  $\sharp$-property iff $B$ has.
    \item [(ii)] If  $A$ has  $\sharp$-property, then it is $\mu$-invariant.
\end{enumerate}
\end{fact}
By Lemma 3.9 and Lemma 3.10 of  \cite{Yao-Zhang-min-flow-II}, we have
\begin{fact}\label{fact-commu-mu-inv}
Any commutative dfg group over $\Qp$ has $\sharp$-property (so is $\mu$-invariant).
\end{fact}

\begin{lemma}\label{lemma-semi-product-sharp}
Let $A$ and $B$ be groups definable over $\Q_p$ with $\sharp$-property.  Then $G=A\rtimes B$  has $\sharp$-property.
\end{lemma}
\begin{proof}
We identify an element $g\in G$ with a pair $(g_A,g_B)$ for $g_A\in A$ and $g_B\in B$. Consider the $\Qp$-definable map  $\eta:B\times A\to A, \ (y,x)\mapsto x^y$. Firstly, we have that pairs  $(x, y)$ and  $(\eta(y^{-1}, x),y)$ are interdefinable over $\Qp$ for $x\in A$ and $y\in B$.

Let $N^*\succ N$ and $p=\tp(a, b/N^*)\in S_G(N^*)$ be $G^{0}(N)$-invariant. We now show that
$p$ is $\mu$-invariant. It is easy to see that $\mu_G(\M)=\mu_A(\M)\rtimes \mu_B(\M)$ and $G^{0}=A^{0}\rtimes B^{0}$.

For every $b'\in B^{0}(N)$,
\[
(\id_A,b')\cdot \tp(a, b/N^*)=\tp(\eta(b',a),b'b/N^*)=\tp(a, b/N^*).
\]
We see from the  correspondence $(x, y)\mapsto (\eta(y^{-1}, x),y)$ that
\[
\tp(\eta(b^{-1}, a),b/N^*)=\tp(\eta((b'b)^{-1},\eta(b', a)),b'b/N^*)=\tp(\eta(b^{-1}, a),b'b/N^*).
\]
Hence, $\tp(b/N^*, \eta(b^{-1}, a))$ is $B^{0}(N)$-invariant. Take any $\delta\in\mu_B(\M)$. By the $\sharp$-property of $B$,  we have that
\[
\tp(\delta b/N^*, \eta(b^{-1}, a))=\tp(b/N^*, \eta(b^{-1}, a)),
\]
i.e.
\[
\tp(\eta(b^{-1}, a) ,\delta b/N^* )=\tp(\eta(b^{-1}, a),b/N^*).
\]
Applying the  correspondence $(x, y)\mapsto (\eta(y^{-1}, x),y)$ again, we have that
\[
\tp(\eta(\delta, a), \delta b/N^*)=\tp(a,b/N^*)=p.
\]

As $\tp(\eta(\delta, a), \delta b/N^*)=p$ is $G^{0}(N)$-invariant,  $\tp(\eta(\delta, a)/N^*,\delta b)$ is $A^{0}(N)$-invariant. Take any $\epsilon\in \mu_A(\M)$. Since $A$ has $\sharp$-property,
$\tp(\epsilon\eta(\delta, a)/N^*,\delta h)=\tp(\eta(\delta, a)/N^*,\delta h)$. Thus,
\[
(\epsilon,\delta)\cdot \tp(a,b/N^*)=\tp(\epsilon\eta(\delta, a), \delta b/N^*)=\tp(\eta(\delta, a), \delta b/N^*)=\tp(a, b/N^*)
\]
for any $(\epsilon,\delta)\in\mu_A(\M)\rtimes\mu_B(\M)$, which completes the proof.
\end{proof}

\begin{proposition}\label{prop-tri-mu-inv}
     $H$ has $\sharp$-property hence is $\mu$-invariant.
\end{proposition}

\begin{proof}
Induction on $\dim(H)$. Now $H=U\rtimes T$ with $U$ unipotent and $T=\Gm^k$.  If $\dim(H)\leq 1$, then $H$ is either $\Ga$ or $\Gm$ or trivial, so $H$ has $\sharp$-property.  If $\dim(U)=0$, then $H=T$ is commutative, and hence, has $\sharp$-property by Fact \ref{fact-commu-mu-inv}. If $\dim(T)>0$, then $\dim(U)<\dim(H)$, so $U$ has  $\sharp$-property by induction hypothesis, and therefore by Lemma \ref{lemma-semi-product-sharp}, $H$ also has $\sharp$-property.  So we assume that $\dim(U)>0$ and $\dim(T)=0$. By Lemma \ref{lemma-uni-semi-pro}, $U=U_1\rtimes \Ga$ with $U_1$ unipotent and $\dim(U_1)=\dim(U)-1$. Applying the induction hypothesis and  Lemma \ref{lemma-semi-product-sharp} again, we see that $U$  has $\sharp$-property. This completes the proof.
\end{proof}

We conclude directly from Fact \ref{fact-dfg-groups}, Fact \ref{fact-virtually-sationary-bad}, Fact \ref{fact-sharp-virtually-mu-inv}, Proposition \ref{prop-stationary iff not bad}, and Proposition \ref{lemma-semi-product-sharp} that

\begin{theorem}\label{thm-stationary-bad-mu-inv}
If $H$ is a dfg group over $\Qp$, then
\begin{enumerate}
    \item [(i)] $H$ is stationary iff $H$ is not bad.
    \item [(ii)] $H$ is $\mu$-invariant.
\end{enumerate}

\end{theorem}


\section{Weakly generics and almost periodics}

In this section $G\sq\M^n$ will be  a group definable over $\Qp$ admitting a $\Qp$-definable short exact sequence
\[
1\rightarrow H\rightarrow G\rightarrow_\pi C\rightarrow 1,
\]
where $C$ is an infinite fsg  group  and $H$ an infinite $\dfg$ group.   Recall that $p$CF has definable skolem functions \cite{Skolem}, so $\pi: G\to C$ has a $\Qp$-definable section $f:C\rightarrow G$, and thus any $g\in G$ can be written uniquely as $f(c)h$  for $c=\pi(g)\in C$ and $h=(f(c))^{-1}g\in H$. Let $\eta: C\times C\rightarrow H$ defined by \[
\eta(c_1,c_2)=f(c_1c_2)^{-1}f(c_1)f(c_2).
\]
Then $\eta$ is also a $\Qp$-definable function, and we have $f(c_1)f(c_2)=f(c_1c_2)\eta(c_1,c_2)$.
For any $N\succ \Qp$ and $\tp(c/N)\in S_C(N)$, by $f(\tp(c/N))$ we mean the type $\tp(f(c)/N)$.

We recall some results in earlier papers first. For the dfg group $H$, we have:
\begin{fact}\label{fact-WG-AP-dfg}\cite{P-Y-dfg groups}
Let $N\succ \Qp$.   Then $\WG(S_H(N))=\AP(S_H(N))$.
\end{fact}
Since $H$ is virtually trigonalizable over $\Qp$, we see from
 Corollary 2.16 of \cite{Yao-tri} that

 \begin{fact}\label{fact-WG-AP-dfg-II}
Let $N\succ \Qp$ and $p\in \WG(S_H(N))$, then every global heir of $p$ is weakly generic.
\end{fact}

For the fsg group $C$, we have:
\begin{fact}\label{fact-WG-AP-fsg}\cite{O-P}
Let $N\succ \Qp$. Then
    \begin{enumerate}
    \item [(i)] $\Gen(S_C(N))$ is not empty.
        \item [(ii)] If $N$ is sufficiently saturated, then  $p\in \Gen(S_C(N))$ iff  every left and right $C(N)$-translate of $p$ is finitely satiafiable in $\Qp$. Namely, $p$ is fsg when $N=\M$.
        \item [(iii)] If $p\in \Gen(S_C(N))$, then $p$ has a unique global generic extension.
      \item [(iv)] If $p\in \Gen(S_C(\Qp))$, then  the  global generic extension of $p$ is precisely the unique global coheir of $p$.
    \end{enumerate}
\end{fact}

\begin{fact}\label{fact-min-flow-restriction}\cite{Sim-VC}
Let $T$ be an NIP theory, $A\prec B\models T$, $X$ an  $A$-definable group defined in the monster model of $T$,  $\pi: S_X(B)\rightarrow S_X(A)$ the canonical restriction map, and $\cal M$  a minimal $X(B)$-subflow of $S_X(B)$. Then $\pi({\cal M})$ is a  minimal $X(A)$-subflow of $S_X(A)$.
\end{fact}

\begin{fact}\label{Semigp-struc}\cite{Newelski-I}
We can equip $S_G(\Qp)$ with a semigroup structure $(S_G(\Qp),*)$: for $p,q\in S_G(\Qp)$, $p*q=\tp(a\cdot b/\Qp)$ with $a$ realizing $ p$ and $b$ realizing the unique heir of $q$ over $\dcl(\Qp,a)$.
\end{fact}

\begin{fact}\cite{Yao-Zhang-min-flow-I}\label{fact-AP=q-p-r}
Let $r\in S_G(\Qp)$. Then $r$ is almost periodic iff $r=f(q)*p*r$ for some  $q\in \Gen(S_C(\Qp))$ and  $p\in \WG(S_H(\Qp))$.
\end{fact}

We can construct a strongly $f$-generic type on $G$ as following:
\begin{fact}\label{fact-construct-strongly-f-generic}\cite{Yao-Zhang-min-flow-I}
 Let  $\tp(c^*/\M )$ be a generic type on $C$ and $\tp(h^*/  \M, c^* )$ is a weakly type on $H$,  strongly $f$-generic over $N\prec\M$. Then $\tp(f(c^*)h^*/\M)$ is a strongly $f$-generic type on $G$ over $N$.
\end{fact}

\begin{proposition}\label{lemma-AP-neq-WG-local}
Suppose that $H$ is bad, then $\WG(S_G(\Qp))\neq \AP(S_G(\Qp))$.

\end{proposition}
\begin{proof}
Let $N$ be a sufficiently saturated (small)  extension of $\Qp$. Let  $\tp(c^*/N )$ be a generic type on $C$. Since $H$ is bad, there is  $h^*\in H$ such that $\tp(h^*/N,c^*)$ is a weakly generic  type on $H$, strongly $f$-generic over $\Qp$, and there is an $\Qp$-definable function $\theta$ such that $\tp(\theta(h^*)/N,c^*)$ is  finitely satisfiable in $\Qp$ but not realized in $\Qp$.

By Fact \ref{fact-construct-strongly-f-generic}, $\tp(f(c^*)h^*/N)$ is a weakly generic  type on $G$, strongly $f$-generic over $\Q_p$. Let $r=\tp(f(c^*)h^*/\Qp)$, then $r\in \WG(S_G(\Qp))$. Suppose that $r\in \AP(S_G(\Qp))$, then by Fact \ref{fact-AP=q-p-r}, there are $q\in\ii$ and $p\in \jj$ such that $r=f(q)*p*r$.

Let $h\in H(N)$ realize $p$. Let $N_1 $ be an $|N|^+$-saturated extension of $N$, and $f(c')h'\in G(N_1)$ realize the unique heir of $r$ over $N$. Since $\tp(h',c'/\Qp)=\tp(h^*,c^*/\Qp)$, we see that $\tp(h'/\Qp,c')$ is weakly generic and $\tp(\theta(h')/\Qp,c')$ is  finitely satisfiable in $\Qp$ but not realized in $\Qp$.  Note that  $\tp(h'/N,c')$ is an heir of $\tp(h'/\Qp,c')$, so by Fact \ref{fact-WG-AP-dfg}   $\tp(h'/N,c')$ is weakly generic.
 Let $c\in C$ realize the unique coheir of $q$ over $N_1$. Then
\[
f(q)*p*r=\tp(f(c)hf(c')h'/M)=\tp(f(cc')\eta(c,c')h^{f(c')}h'/M).
\]
Now we have the following:
\begin{description}

\item [(*)] $\tp(ah'/N)$ is  definable over $\Q_p$ for each $a\in H(N)$. This is because $\tp(h'/N)$ is the unique heir of $\tp(h'/\Qp)$, so it is definable over $\Q_p$.  By Fact  \ref{fact-strongly-f-generic} and Fact \ref{fact-WG-AP-dfg-II} , it is strongly $f$-generic over $\Q_p$. Take any $a\in H(N)$. then by Remark \ref{rmk-strongly-f-generic-inv}  $\tp(ah'/N)$ is $\Q_p$-invariant. Since $\tp(ah'/N)$ is definable (over $N$), it is definable over $\Q_p$.

\item [(**)] $\tp(h'/N,c,c')$ is an heir of $\tp(h'/\Qp,c')$. By Fact \ref{fact-WG-AP-dfg-II}, $\tp(h'/N,c,c')$ is weakly generic so is $H^{0}(\dcl(N,c,c'))$-invariant.
\end{description}
Since $N$ is sufficiently saturated, $H(N)$ meets every coset of $H^0$. Take $a_0\in H(N)$ such that
\[
a_0/H^0=\eta(c,c')h^{f(c')}h/H^0.
\]
By (**), we have
\[
\tp( \eta(c,c')h^{f(c')}hh'/N,c,c')=\tp(a_0h'/N,c,c').
\]
By (*), we have
\[
\tp( \eta(c,c')h^{f(c')}hh'/N )=\tp(a_0h'/N).
\]
is definable over $\Q_p$.

On the other side,  $\tp(cc'/N_1)$ is a generic type
on $C$, so $\tp(cc'/N)$ is finitely satisfiable in $\Qp$. By Fact \ref{fact-distal}, we see that  $\tp(cc'/N,\eta(c,c')h^{f(c')}hh')$ is finitely satisfiable in $\Qp$. Since
\[
r=\tp(f(c')\cdot h'/\Qp)=\tp\big(f(cc')\cdot\big(\eta(c,c')h^{f(c')}hh'\big)/\Qp\big),
\]
we conclude that  $\tp(c'/\Qp,\theta(h'))$ is finitely satisfiable in $\Qp$.  Our assumption says that $\tp(\theta(h')/\Qp, c')$ is also finitely satisfiable in $\Qp$, which contradicts to Fact \ref{fact-heir-neq-coheir}.
\end{proof}

\begin{corollary}\label{lemma-AP-neq-WG-global}
Suppose that $H$ is bad, then $\WG(S_G(\M))\neq \AP(S_G(\M))$.
\end{corollary}
\begin{proof}
 Let $\tp(c^*/\M)$ be a generic type on $C$. Since $H$ is bad, there is $h^*\in H$ such that $\tp(h^*/\M,c^*)$ is a weakly generic  type on $H$, strongly $f$-generic over $\Qp$. Then $\bar r=\tp(f(c^*)h^*/\M)$ is weakly generic type on $G$. If $r$ is almost periodic, then by Fact \ref{fact-min-flow-restriction}, $r=\tp(f(c^*)h^*/\Qp)$ is a almost periodic type over $\Qp$.  But Lemma \ref{lemma-AP-neq-WG-local} shows that $r$ is not almost periodic. A contradiction.
\end{proof}

Recall the main result of \cite{Yao-Zhang-min-flow-I} is:
\begin{fact}\label{fact-sta-AP=WG}
\begin{enumerate}
\item [(i)] $G$ is stationary iff $H$ is stationary.
\item [(ii)] $G$ is stationary iff $G$ has boundedly may global weakly generic types.
\item [(ii)] (Local case)  If $G$ be stationary, then
\[
\WG(S_G(\Qp ))=\AP(S_G(\Qp ))=f(\AP(S_C(\Qp )))*\AP(S_H(\Qp)).
\]
\item [(iii)] (Global case)  If $G$ be stationary, then $\AP(S_G(\M))=\WG(S_G(\M))$.
\end{enumerate}
\end{fact}

Summarizing Theorem \ref{thm-stationary-bad-mu-inv}, Proposition \ref{lemma-AP-neq-WG-local}, Corollary \ref{lemma-AP-neq-WG-global}, and Fact \ref{fact-sta-AP=WG}, we have:
\begin{theorem}
Let $G$ be a group definable over $\Qp$ which  admits a $\Q_p$-definable short exact sequence
\[
1\rightarrow H\rightarrow G\rightarrow_\pi C\rightarrow 1,
\]
where $C$ is an infinite fsg  group  and $H$ an infinite $\dfg$ group. Then the following are equivalent:
\begin{enumerate}
\item [(i)] $H$ is stationary;
\item [(ii)] $H$ is not bad;
\item [(iii)] $G$ is stationary;
\item [(iv)] $G$ has boundedly many global weakly generic types;  
\item [(v)] $\WG(S_G(\Qp))=\AP(S_G(\Qp))$;
\item [(vi)] $\WG(S_G(\M))=\AP(S_G(\M))$.
\end{enumerate}
\end{theorem}

\subsubsection*{Acknowledgments.}

The research is supported by The National Social Science Fund of China (Grant No. 20CZX050). We would like to thank Will Johnson who provided the proof of Lemma \ref{lemma-uni-semi-pro}.

\end{document}